\documentclass[11pt]{article}
\usepackage{amssymb,amsmath}

\newcommand{\be}{\begin{equation}}
\newcommand{\ee}{\end{equation}}
\newcommand{\bea}{\begin{eqnarray}}
\newcommand{\eea}{\end{eqnarray}}

\newtheorem{corollary}{Corollary}
\newtheorem{remark}{Remark}
\newtheorem{example}{Example}
\newtheorem{lemma}{Lemma}

\def\1#1{^{(#1)}}

\begin{document}
\title{On Numerical Experiments with Symmetric Semigroups Generated by 
Three Elements and Their Generalization}
\author{Leonid G. Fel\\  
\\Department of Civil Engineering, Technion, Haifa 3200, Israel}
\vspace{-.2cm}
\date{\today}   
\maketitle
\def\be{\begin{equation}}
\def\ee{\end{equation}}
\def\bea{\begin{eqnarray}}
\def\eea{\end{eqnarray}}
\def\p{\prime}
\vspace{-.5cm}
\begin{abstract}
We give a simple explanation of numerical experiments of V. Arnold with two
sequences of symmetric numerical semigroups, ${\sf S}(4,6+4k,87-4k)$ and ${\sf 
S}(9,3+9k,85-9k)$ generated by three elements. We present a generalization of 
these sequences by numerical semigroups ${\sf S}(r_1^2,r_1r_2+r_1^2k,r_3-r_1^2
k)$, $k\in{\mathbb Z}$, $r_1,r_2,r_3\in{\mathbb Z}^+$, $r_1\geq 2$ and $\gcd(
r_1,r_2)=\gcd(r_1,r_3)=1$, and calculate their universal Frobenius number $\Phi
(r_1,r_2,r_3)$ for the wide range of $k$ providing semigroups be symmetric.
We show that this kind of semigroups admit also nonsymmetric representatives. 
We describe the reduction of the minimal generating sets of these semigroups up 
to $\{r_1^2,r_3-r_1^2k\}$ for sporadic values of $k$ and find these values by 
solving the quadratic Diophantine equation.
\end{abstract} 
%=================================================================
\section{Introduction}\label{s1}
In experiments with Frobenius numbers $F(d_1,d_2,d_3)$ of numerical semigroups, 
generated by a tuple of three elements $\{d_1,d_2,d_3\}$, V. Arnold has 
mentioned two strange arithmetic facts 
\footnote{Throughout the paper we use the term {\em Frobenius number} whose 
standard definition dates back to G. Frobenius, I. Schur and A. Brauer 
\cite{brau42} and denotes the largest integer that is not representable as a 
linear combination with nonnegative integer coefficients of a given tuple of 
positive integers $\{d_1,\ldots, d_m\}$, $\gcd(d_1,\ldots,d_m)=1$. V. Arnold 
\cite{arn1} had used a different definition of this term, so in (\ref{a1}) and 
(\ref{a2}) he got numbers 90 and 168 instead of 89 and 167.}
(see \cite{arn1}, Remark 1), 
\begin{eqnarray}
F(4,6+4k,87-4k)&=&89\;,\;\;\;\;k=0,1,\ldots ,14,\;k\neq 8\;,\label{a1}\\
F(9,3+9k,85-9k)&=&167\;,\;\;k=1,\ldots ,7\;.\label{a2}
\end{eqnarray}
In sections \ref{s2} and \ref{s3} we give a simple proof of these statements. 
In fact, we prove
\begin{eqnarray}
F(4,6+4k,87-4k)=89\;,\;\;-1\leq k\leq 14\;;\;\;\;
F(9,3+9k,85-9k)=167\;,\;\;0\leq k\leq 7\;.\label{a3}
\end{eqnarray}
The proof is based on observation that two sequences of triples, 
\begin{eqnarray}
\{4,6+4k,87-4k\},\;0\leq k\leq 14\;,\;\;\;\mbox{and}\;\;\;
\{9,3+9k,85-9k\},\;1\leq k\leq 6\;,\label{a4}
\end{eqnarray}
present the generators of symmetric numerical semigroups generated by three 
elements. The case $k=-1$ in the 1st triple and the cases $k=0$, $k=7$ in the 
2nd triple are special and reduce the semigroups, which are generated by three 
elements, up to symmetric semigroups, which are generated by two elements. In 
sections \ref{s4} and \ref{s5} we generalize both examples (\ref{a4}) to most 
generic triple and analyze its associated symmetric semigroups. In section 
\ref{s6} we discuss a case when this generating triple is reduced up to 
generating pair, and values of their elements are coming by finding the 
integer points in plane algebraic curve of degree 2.
%=================================================================
\subsection{Basic Facts on Numerical Semigroups ${\sf S}\left({\bf d}^3\right)$}
\label{s11}
Following \cite{fe1} we recall basic definitions and known facts on algebra of 
the numerical semigroups generated by $m$ elements which are necessary here, 
and focus on their symmetric subsets. For short we denote the generating tuple 
$(d_1,\ldots,d_m)$ by ${\bf d}^m$.

A semigroup ${\sf S}\left({\bf d}^m\right)=\left\{s\in{\mathbb Z}^+\cup\{0\}\;|
\;s=\sum_{i=1}^m x_i d_i,\;x_i\in {\mathbb Z}^+\cup\{0\}\right\}$ is said to be 
generated by minimal set of $m$ natural numbers, 
\begin{eqnarray}
\gcd(d_1,\ldots,d_m)=1\;,\;\;\;\;\min(d_1,\ldots,d_m)\geq m\;,\label{a0}
\end{eqnarray}
if neither of its elements is linearly representable by the rest of elements. 
Throughout the paper we call such semigroups $m$-dimensional (mD). Denote by 
$\Delta\left({\bf d}^m\right)$ the complement of ${\sf S}\left({\bf d}^m
\right)$ in ${\mathbb Z}^+$, i.e. $\Delta\left({\bf d}^m\right)={\mathbb Z}^+
\setminus{\sf S}\left({\bf d}^m\right)$ and call it a set of gaps. The Frobenius
number of semigroup ${\sf S}\left({\bf d}^m\right)$ is defined as follows, 
\begin{eqnarray}
F\left({\bf d}^m\right)=\max\Delta\left({\bf d}^m\right)\;.\label{a5}
\end{eqnarray}
A semigroup ${\sf S}\left({\bf d}^m\right)$ is called symmetric if for any 
integer $s$ the following condition holds: iff $s\in {\sf S}\left({\bf d}^m
\right)$ then $F\left({\bf d}^m\right)-s\not\in {\sf S}\left({\bf d}^m\right)$. 
Otherwise ${\sf S}\left({\bf d}^m\right)$ is called nonsymmetric. Notably that 
all semigroups ${\sf S}\left(d_1,d_2\right)$, $\min(d_1,d_2)\geq 2$, generated 
by two elements, i.e. two-dimensional (2D) semigroups, are symmetric. Combining 
the last fact with the early statement of Watanabe on the symmetric semigroups 
of dimension $m\geq 3$ (see \cite{wa73}, Lemma 1) we come to important 
statement.
\begin{lemma}\label{lem1}
Let ${\sf S}\left(c_1,c_2\right)$ be a numerical semigroup, $a$ and $b$ be
positive integers, $\gcd(a,b)=1$. If $a\in {\sf S}\left(c_1,c_2\right)$, then 
the semigroup ${\sf S}\left(bc_1,bc_2,a\right)$ is symmetric.
\end{lemma}
\begin{remark}\label{rem1}
Note that a requirement $a\in {\sf S}\left(c_1,c_2\right)$ can be provided in 
two ways. First, it is satisfied if the Frobenius number of 2D semigroup 
${\sf S}(c_1,c_2)$ is exceeded by the third generator '$a$',
\begin{eqnarray}
a\geq C(c_1,c_2)=1+F(c_1,c_2)=(c_1-1)(c_2-1)\;,\label{a20}
\end{eqnarray}
where $C(c_1,c_2)$ denotes a conductor of semigroup ${\sf S}(c_1,c_2)$. The last
equality in (\ref{a20}) comes due to the known Sylvester formula \cite{herz70}. 
In the case $a<(c_1-1)(c_2-1)$ there is another way to provide the containment 
$a\in {\sf S}(c_1,c_2)$, namely, to have the number '$a$' among the nongaps of 
semigroup ${\sf S}(c_1,c_2)$, i.e. $a\in {\sf S}(c_1,c_2)\cap [0,F(c_1,c_2)]$. 
The last requirement is much harder to verify than (\ref{a20}).
\end{remark}
The last case will be also observed in further calculation (see section 
\ref{s3}, a case $k=6$ for the triple $\{9,3+9k,85-9k\}$). A powerful tool to 
study the symmetric numerical semigroups is the Herzog formula \cite{herz70} 
for the Frobenius number (for details see section 6.2 in \cite{fe1}). Being 
adapted for symmetric semigroup ${\sf S}\left(bc_1,bc_2,a \right)$ in Lemma 
\ref{lem1}, it looks as follows,
\begin{eqnarray}
F\left(bc_1,bc_2,a\right)=bc_1c_2+ab-(bc_1+bc_2+a)\;.\label{a6}
\end{eqnarray}
Keeping in mind Lemma \ref{lem1} consider in details the two sequences of 
triples given in (\ref{a4}).
%=================================================================
\section{Symmetric Semigroups ${\sf S}(4,6+4k,87-4k)$}\label{s2}
The triple $\{4,6+4k,87-4k\}$ has always two relative non-prime generators 
$bc_1$ and $bc_2$, namely, $c_1=2$, $c_2=3+2k$ and $b=2$. In order to satisfy 
Lemma \ref{lem1} we have to provide the containment $87-4k\in {\sf S}(2,3+2k)$. 
By (\ref{a0}) the requirement $3+2k\geq 2$ brings us to the lower bound for 
$k$, $k\geq 0$. The upper bound comes by another claim for conductor 
$C(2,3+2k)$ of semigroup ${\sf S}(2,3+2k)$,
\begin{eqnarray}
87-4k\geq C(2,3+2k)=2\cdot (3+2k)-(2+3+2k)+1\;\;\rightarrow\;\;85\geq 6k\;\;
\rightarrow\;\;k\leq 14\;.\label{a7}
\end{eqnarray}
Thus, by Lemma \ref{lem1} the numerical semigroups ${\sf S}(4,6+4k,87-4k)$, 
$0\leq k\leq 14$, are symmetric. Applying (\ref{a6}) we get
\begin{eqnarray}
F(4,6+4k,87-4k)=4\cdot (3+2k)+2\cdot (87-4k)-(4+6+4k+87-4k)=89\;.\label{a8}
\end{eqnarray}

The higher values of $k$ are bounded by the claim (\ref{a0}): $87-4k\geq 3$ 
that gives $k\leq 21$. The corresponding semigroups ${\sf S}(4,6+4k,87-4k)$, 
$15\leq k\leq 21$, are isomorphic to the 2D symmetric semigroups:
\begin{eqnarray}
{\sf S}(4,66,27)={\sf S}(4,27),\;\;{\sf S}(4,70,23)={\sf S}(4,23),\;\;
{\sf S}(4,74,19)={\sf S}(4,19),\;\;{\sf S}(4,78,15)={\sf S}(4,15),\nonumber\\
{\sf S}(4,82,11)={\sf S}(4,11),\;\;{\sf S}(4,86,7)={\sf S}(4,7),\;\;
{\sf S}(4,90,3)={\sf S}(4,3)\;.\nonumber
\end{eqnarray}
Their Frobenius numbers can be found by Sylvester formula,
\begin{eqnarray}
F(4,66,27)=77,\;\;F(4,70,23)=65,\;\;F(4,74,19)=53,\;\;F(4,78,15)=41,\nonumber\\
F(4,82,11)=29,\;\;F(4,86,7)=17,\;\;F(4,90,3)=5\;.\nonumber
\end{eqnarray}
The case $k=-1$ is a special one. It corresponds to semigroup ${\sf S}(4,2,91)$
with non-minimal generating set $\{4,2,91\}$. It can be reduced up to 
$\{2,91\}$ which generates a semigroup ${\sf S}(2,91)$. The Frobenius number 
of the latter semigroup follows by Sylvester formula, $F(2,91)=89$.
%=================================================================
\section{Symmetric Semigroups ${\sf S}(9,3+9k,85-9k)$}\label{s3}
The triple $\{9,3+9k,85-9k\}$ has always two relative non-prime generators 
$bc_1$ and $bc_2$, namely, $c_1=3$, $c_2=3k+1$ and $b=3$. In order to satisfy 
Lemma \ref{lem1} we have to provide the containment $85-9k\in {\sf S}(3,3k+1)$.
By (\ref{a0}) the requirement $3k+1\geq 2$ brings us to the lower bound for $k$,
$k\geq 1$. The upper bound comes by another claim for conductor $C(3,3k+1)$ of 
semigroup ${\sf S}(3,3k+1)$,
\begin{eqnarray}
85-9k\geq C(3,3k+1)=3\cdot (3k+1)-(3+3k+1)+1\;\;\rightarrow\;\;85\geq 15k\;\;
\rightarrow\;\;k\leq 5\;.\label{a9}
\end{eqnarray}
Thus, by Lemma \ref{lem1} the numerical semigroups ${\sf S}(9,3+9k,85-9k)$, 
$1\leq k\leq 5$, are symmetric. Applying (\ref{a6}) we get
\begin{eqnarray}
F(9,3+9k,85-9k)=9\cdot (3k+1)+3\cdot (85-9k)-(9+9k+3+85-9k)=167\;.\label{a10}
\end{eqnarray}
A case $k=6$ gives rise to another symmetric semigroup ${\sf S}(9,57,31)$ 
which satisfies Lemma \ref{lem1}: $31\in {\sf S}(3,19)$, however $31<C(3,19)$.
Making use of (\ref{a10}) we get $F(9,57,31)=167$.

Finally, two other cases $k=0$ and $k=7$ give rise to 3D semigroups ${\sf S}(9,
3,85)$ and ${\sf S}(9,66,22)$ with non-minimal generating sets $\{9,3,85\}$ and 
$\{9,66,22\}$, respectively. However, they can be reduced up to the 2D 
semigroups ${\sf S}(3,85)$ and ${\sf S}(9,22)$, respectively. The Frobenius 
numbers of the two last semigroups follow by Sylvester formula, $F(9,3,85)=167$ 
and $F(9,66,22)=167$.

The higher values of $k$ are bounded by the claim (\ref{a0}): $85-9k\geq 3$ 
that gives $k\leq 9$. The corresponding semigroups ${\sf S}(9,3+9k,85-9k)$, 
$k=8,9$, are isomorphic to the 2D symmetric semigroups:
\begin{eqnarray}
{\sf S}(9,75,13)={\sf S}(9,13)\;,\;\;{\sf S}(9,84,4)={\sf S}(9,4)\;.\nonumber
\end{eqnarray}
Their Frobenius numbers follow  by Sylvester formula, $F(9,75,13)=95$, 
$F(9,84,4)=23$.

It is worth to mention that in the whole range of varying parameter $k$ with 
the values of the triples' elements exceeding 1 both sequences of these 
triples in sections \ref{s2} and \ref{s3} give rise only to symmetric 
semigroups either three-dimensional or two-dimensional,
\begin{eqnarray}
\{4,6+4k,87-4k\},\;-1\leq k\leq 21\;,\;\;\;\mbox{and}\;\;\;
\{9,3+9k,85-9k\},\;0\leq k\leq 9\;.\label{a01}
\end{eqnarray}
This observation is important not less than the claim (\ref{a3}) on 
universality of the Frobenius numbers 89 and 167. However the range of 
application of (\ref{a01}) is much wider than (\ref{a4}). 
%=================================================================
\section{Numerical Semigroups ${\sf S}(r_1^2,r_1r_2+r_1^2k,r_3-r_1^2k)$}
\label{s4}
In this section we generalize both examples discussed by V. Arnold in 
\cite{arn1}. For the first glance, a most generic triple is of the form,
\begin{eqnarray}
\left\{u^2v^2,u^2vw+u^2v^2k,t-u^2v^2k\right\},\left\{\begin{array}{l}k\in
{\mathbb Z},\\u,v,w,t\in{\mathbb Z}^+\end{array}\right.,\left\{\begin{array}
{l}\gcd(u,v)=\gcd(u,w)=\gcd(v,w)=1,\\\gcd(u,t)=\gcd(v,t)=1,\;\;\;uv\geq 2.
\end{array}\right.\nonumber
\end{eqnarray}
However, by comparison with Arnold's examples, the last triple has one serious 
lack. Indeed,  consider a symmetric semigroup  ${\sf S}\left(u^2v^2,u^2vw+
u^2v^2k,t-u^2v^2k\right)$ and calculate by formula (\ref{a6}) its Frobenius 
number,
\begin{eqnarray}
F(u,v,w,t,k)=(t+u^2vw)(v-1)-u^2v^2+ku^2v^3(1-u^2)\;.\label{a50a}
\end{eqnarray}
In contrast to examples in \cite{arn1}, an expression in (\ref{a50a}) is 
dependent on $k$. This dependence disappears iff $u=1$. The generating triples 
of only such kind will be a subject of interest in this article. Henceforth,
denote $v=r_1$, $w=r_2$, $t=r_3$  and consider a triple which is governed by 
three parameters, $r_1$, $r_2$ and $r_3$,
\begin{eqnarray}
\left\{r_1^2,r_1r_2+r_1^2k,r_3-r_1^2k\right\},\;\;k\in{\mathbb Z},\;r_1,r_2,r_3
\in{\mathbb Z}^+,\;r_1\geq 2\;\;\mbox{and}\;\;\gcd(r_1,r_2)=\gcd(r_1,r_3)=1.\;
\label{a50}
\end{eqnarray}
In new notations $r_1$, $r_2$ and $r_3$ and by $u=1$ formula (\ref{a50a}) reads
\begin{eqnarray}
\Phi\left(r_1,r_2,r_3\right)=(r_1-1)(r_1r_2+r_3)-r_1^2\;.\label{a23}
\end{eqnarray}
There are two different ways to symmetrize the 3D numerical semigroup ${\sf S}
(r_1^2,r_1r_2+r_1^2k,r_3-r_1^2k)$. The 1st way is to choose $k$ such that the 
necessary conditions in Lemma \ref{lem1} be satisfied. The 2nd way is to choose 
$k$ such that the generating triple is non-minimal, i.e. one of its elements is
linearly representable by the rest of elements. In other words, one can arrive
at symmetric semigroup preserving the dimension 3 of generic semigroup or
reducing it by 1. 

Unfortunately, a complete analysis of symmetrization of numerical semigroup 
generated by the triple (\ref{a50}) encounter a serious difficulty in both ways 
of its performing. This is related to non-analytic nature of both containments 
$r_3-r_1^2\bar{k}\in{\sf S}(r_1,r_2+r_1\bar{k})$ and $r_1r_2+r_1^2\tilde{k}\in
{\sf S}(r_1^2,r_3-r_1^2\tilde{k})$. In other words, one cannot write the 
explicit formulas of $\bar{k}$ and $\tilde{k}$ via $r_1,r_2,r_3$ for the whole 
set of nongaps for both semigroups ${\sf S}(r_1,r_2+r_1\bar{k})$ and ${\sf S}(
r_1^2,r_3-r_1^2\tilde{k})$. What we can do only to make use of (\ref{a20}) 
providing the ranges of $\bar{k}$ and $\tilde{k}$ when the elements $r_3-r_1^2
\bar{k}$ and $r_1r_2+r_1^2\tilde{k}$ exceed the Frobenius numbers of 
corresponding semigroups, respectively. According to Remark \ref{rem1} this 
symmetrizes an initial semigroup ${\sf S}(r_1^2,r_1r_2+r_1^2k,r_3-r_1^2k)$ at
$k=\bar{k},\tilde{k}$.

In section \ref{s41} we find the range of $k$-values for the sequence of 
symmetric semigroups generated by the triple (\ref{a50}) with equal Frobenius 
numbers (\ref{a23}). In section \ref{s42} we give an affirmative answer to 
another question: whether the sequence (\ref{a50}) does contain also a triple 
associated with nonsymmetric semigroups.
%=================================================================
\subsection{Symmetric semigroups and special values of $k$}\label{s41}
We start with the 1st way of symmetrization and assume that $r_1k+r_2\neq 1$. 
By Lemma \ref{lem1} a numerical semigroup ${\sf S}(r_1^2,r_1r_2+r_1^2k,
r_3-r_1^2k)$ is symmetric if the following containment holds, $r_3-r_1^2k\in
{\sf S}(r_1,r_1k+r_2)$. By (\ref{a0}) it brings us necessarily to the two 
inequalities imposed onto generators,
\begin{eqnarray}
r_1k+r_2\geq 2\;,\;\;\;\;r_3-r_1^2k\geq 3\;.\label{a21}
\end{eqnarray}
Denote two special values of $k$,
\begin{eqnarray}
k_1=\frac{2-r_2}{r_1}\;,\;\;\;\;k_2=\frac{r_3-3}{r_1^2}\;,\label{a21b}
\end{eqnarray}
and find a range of $k$ where both inequalities (\ref{a21}) do not contradict 
each other,
\begin{eqnarray}
\mbox{if}\;\;\;\;\;k_1\leq k_2\;\;\;\;\;\mbox{then}\;\;\;\;\;k_1\leq k\leq k_2
\;.\label{a31}
\end{eqnarray}
On the other hand, 
\begin{eqnarray}
\mbox{if}\;\;\;\;\;k_1>k_2\;\;\;\;\;\mbox{or}\;\;\;\;\;k\leq k_1\;\;\;\;\;
\mbox{or}\;\;\;\;\;k\geq k_2\;,\label{a31a}
\end{eqnarray}
then the corresponding $k$ does not provide the necessary requirement 
(\ref{a0}).

Apply Lemma \ref{lem1} to semigroup ${\sf S}(r_1^2,r_1r_2+r_1^2k,r_3-r_1^2k)$.
Being 2-dimensional, the symmetric semigroup ${\sf S}(r_1,r_1k+r_2)$ is 
associated with Frobenius number according to Sylvester formula. Following 
(\ref{a20}) write an inequality 
\begin{eqnarray}
r_3-r_1^2k\geq 1+F(r_1,r_1k+r_2)=(r_1-1)(r_1k+r_2-1)\;.\label{a22a}
\end{eqnarray}
It gives rise to another special value of $k$,
\begin{eqnarray}
k\leq k_3\;,\;\;\;\;\;k_3=\frac{r_3-(r_1-1)(r_2-1)}{r_1(2r_1-1)}\;,\label{a22}
\end{eqnarray}
where our 3D semigroup is symmetric. If the inequality (\ref{a22a}) is broken,
\begin{eqnarray}
r_3-r_1^2k\leq F(r_1,r_1k+r_2)\;,\;\;\;\;&\mbox{or}&\;\;\;\;k\geq k_3+\frac1{
r_1(2r_1-1)}\;,\label{a22q}
\end{eqnarray}
then a containment $r_3-r_1^2k\in {\sf S}(r_1,r_1k+r_2)$ can be still provided 
if $r_3-r_1^2k$ is a nongap of semigroup ${\sf S}(r_1,r_1k+r_2)$. Note that 
inequality (\ref{a22q}) admits also the existence of nonsymmetric semigroups
generated by triple (\ref{a50}) if $r_3-r_1^2k$ is a gap of 
${\sf S}(r_1,r_1k+r_2)$.

Let us find the common range of $k$ which is consistent with (\ref{a31}), 
(\ref{a31a}), (\ref{a22}) and (\ref{a22q}) and dependent on interrelationships 
between $k_1$, $k_2$ and $k_3$. By comparison of expressions (\ref{a21b}) and 
(\ref{a22}) for $k_1$, $k_2$ and $k_3$ we find the constraints when these 
relationships are valid. Below we list these relationships presented in terms 
of $r_1$, $r_2$ and $r_3$.
\begin{eqnarray}
k_1\leq k_2\leq k_3\;,\;\;\;\;&\mbox{or}&\;\;\;\;2r_1+3\leq r_1r_2+r_3\leq
\frac{r_1^2+5r_1-3}{r_1-1}\;,\label{a22b}\\
k_1\leq k_3\leq k_2\;,\;\;\;\;&\mbox{or}&\;\;\;\;\frac{r_1^2+5r_1-3}{r_1-1}\leq
r_1r_2+r_3\;,\;\;\;\;3r_1-1\leq r_1r_2+r_3\;,\label{a22c}\\
k_3\leq k_1\leq k_2\;,\;\;\;\;&\mbox{or}&\;\;\;\;2r_1+3\leq r_1r_2+r_3\leq 3r_1
-1\;.\label{a22e}
\end{eqnarray}
%=================================================================
\subsubsection{Semigroup's reduction: ${\sf S}(r_1^2,r_1r_2+r_1^2k,r_3-r_1^2k)
\;\rightarrow\;{\sf S}(r_1,r_3-r_1^2k)$}\label{s411}
Consider the case $r_1k+r_2=1$. Indeed, by this relation the two first 
generators of the triple (\ref{a50}) become linearly dependent, and therefore 
the 3D numerical semigroup is reduced up to the 2D semigroup ${\sf S}(r_1,r_3-
r_1^2k)$ which is always symmetric. Summarizing these requirements we conclude 
that a numerical semigroup ${\sf S}(r_1^2,r_1r_2+r_1^2k_4,r_3-r_1^2k_4)$ is 
symmetric if 
\begin{eqnarray}
k_4=\frac{1-r_2}{r_1}\in {\mathbb Z}\;.\label{a51}
\end{eqnarray}
The corresponding generator $r_3-r_1^2k_4$ and the Frobenius number 
$F(r_1,r_3-r_1^2k_4)$ read
\begin{eqnarray}
r_3-r_1^2k_4=r_3+r_1r_2-r_1\;,\;\;\;F(r_1,r_3-r_1^2k_4)=(r_3+r_1r_2)(r_1-1)-
r_1^2\;.\label{a27b}
\end{eqnarray}
Note that $k_4=k_1-1/r_1$, i.e. $k_1-k_4\leq 1/2$. In fact, this expands the 
range (\ref{a31}) of existence of symmetric numerical semigroups generated by 
the triple (\ref{a50}) up to $k_4\leq k\leq k_2$. Note that two Frobenius 
numbers $F(r_1,r_3-r_1^2k_4)$ and $\Phi(r_1,r_2,r_3)$ given by (\ref{a27b}) 
and (\ref{a23}) coincide. If $r_2=1$ then there always exists 2D semigroup 
${\sf S}(r_1,r_3)$ which comes by putting $k=0$ into (\ref{a50}).
%=================================================================
\subsubsection{Semigroup's reduction: ${\sf S}(r_1^2,r_1r_2+r_1^2k,r_3-r_1^2k)
\;\rightarrow\;{\sf S}(r_1^2,r_3-r_1^2k)$}\label{s412}
Next, consider the case $r_1r_2+r_1^2k\in {\sf S}(r_1^2,r_3-r_1^2k)$ and find
the value $k_5$ such that for all $k>k_5$ the above containment is provided.
For this purpose, in accordance with (\ref{a20}), we have to satisfy the 
following inequality,
\begin{eqnarray}
r_1r_2+r_1^2k>F(r_1^2,r_3-r_1^2k)=(r_3-r_1^2k)(r_1^2-1)-r_1^2\;.\label{a53}
\end{eqnarray}
It gives another special value of $k$,
\begin{eqnarray}
k>k_5\;,\;\;\;\;\;k_5=\frac{(r_3-1)r_1^2-(r_3+r_1r_2)}{r_1^4}\;.\label{a54}
\end{eqnarray}
By (\ref{a21b}) and (\ref{a54}) it follows 
\begin{eqnarray}
k_2-k_5=\frac{r_1r_2+r_3-2r_1^2}{r_1^4}\;,\;\;\;\;\;\;\mbox{i.e.}\;\;\;\;
k_2\geq k_5\;\;\;\;\mbox{iff}\;\;\;\;r_1r_2+r_3\geq 2r_1^2\;.\label{a54a}
\end{eqnarray}
Find a value $k_6$ where the Frobenius number $F(r_1^2,r_3-r_1^2k_6)$ 
coincides with $\Phi(r_1,r_2,r_3)$,
\begin{eqnarray}
F(r_1^2,r_3-r_1^2k_6)=(r_1-1)(r_1r_2+r_3)-r_1^2\;\;\;\rightarrow\;\;\;
k_6=\frac{r_3-r_2}{r_1(r_1+1)}\;.\label{a55}
\end{eqnarray}
%=================================================================
\subsubsection{Semigroup's reduction: ${\sf S}(r_1^2,r_1r_2+r_1^2k,r_3-r_1^2k)
\;\rightarrow\;{\sf S}(r_1r_2+r_1^2k,r_3-r_1^2k)$}\label{s413}
Finally, consider the case $r_1^2\in {\sf S}(r_1r_2+r_1^2k,r_3-r_1^2k)$ and 
find the $k$-values such that the above containment is provided. In 
accordance with (\ref{a20}), we have to satisfy the following inequality,
\begin{eqnarray}
r_1^2\geq F(r_1r_2+r_1^2k,r_3-r_1^2k)+1\;.\label{a55a}
\end{eqnarray}
It gives two other special values of $k$,
\begin{eqnarray}
k\leq k_7\;\;\;\;\mbox{or}\;\;\;\;k_8\leq k\;,\;\;\;\mbox{where}\;\;\;\;\;\;
k_8-k_7=\frac{\sqrt{(r_3+r_1r_2-2)^2-4r_1^2}}{r_1^2}\;,\label{a55b}
\end{eqnarray}
\begin{eqnarray}  
k_7=\frac{r_3-r_1r_2-\sqrt{(r_3+r_1r_2-2)^2-4r_1^2}}{2r_1^2}\;,\;\;\;\;
k_8=\frac{r_3-r_1r_2+\sqrt{(r_3+r_1r_2-2)^2-4r_1^2}}{2r_1^2}\;.\label{a55c}
\end{eqnarray}
By (\ref{a55b}) we have $k_8\geq k_7$, if $r_3+r_1r_2\geq 2+2r_1$, otherwise an 
inequality (\ref{a55a}) holds for any $k$. Making use of formulas (\ref{a21b}), 
(\ref{a55c}) and calculating two differences, $k_8-k_2$ and $k_1-k_7$, we get
\begin{eqnarray}
a)\;\;k_8\geq k_2\;\;\;\;\mbox{if}\;\;\;\;r_3+r_1r_2\geq 4+\frac{r_1^2}{2}\;;
\;\;\;\;\;\;\;b)\;\;k_1\geq k_7\;\;\;\;\mbox{if}\;\;\;\;r_3+r_1r_2\geq
\frac{5r_1^2-1}{2r_1-1}\;.\label{a55d}
\end{eqnarray}
By comparison of criteria in (\ref{a54a}) and (\ref{a55d}) we obtain
\begin{eqnarray}
\mbox{if}\;\;\;k_2\geq k_5\;,\;\;\;\mbox{then}\;\;\;k_7\leq k_1\;\;\;
\mbox{and}\;\;\;k_2\leq k_8\;.\label{a55e}
\end{eqnarray}
In other words, if $r_1r_2+r_3\geq 2r_1^2$ then the semigroup's reduction 
${\sf S}(r_1^2,r_1r_2+r_1^2k,r_3-r_1^2k)\rightarrow{\sf S}(r_1r_2+r_1^2k,
r_3-r_1^2k)$ cannot be observed. On the other hand,
\begin{eqnarray}
\mbox{if}\;\;\;k_1\leq k_7\;\;\;\mbox{and}\;\;\;k_8\leq k_2\;,\;\;\;
\mbox{then}\;\;\;k_2\leq k_5\;.\label{a55f}
\end{eqnarray}
However, the opposite relationship is not always true,
\begin{eqnarray}
\mbox{if}\;\;\;\max\left\{4+\frac{r_1^2}{2}\;,\;\frac{5r_1^2-1}{2r_1-1}\right\}
\leq r_1r_2+r_3\leq 2r_1^2\;,\;\;\;\;\mbox{then}\;\;\;\;\left\{\begin{array}{l}
k_7\leq k_1,\;k_2\leq k_8\;,\\k_2\leq k_5\;.\end{array}\right.\label{a55h}
\end{eqnarray}
Find two other values $k_9$ and $k_{10}$ where the Frobenius numbers $F(r_1r_2
+r_1^2k_9,r_3-r_1^2k_9)$ and $F(r_1r_2+r_1^2k_{10},r_3-r_1^2k_{10})$ coincide 
with $\Phi(r_1,r_2,r_3)$,
\begin{eqnarray}   
k_9=\frac{r_3-r_1}{r_1^2}\;,\;\;\;\;k_{10}=\frac{1-r_2}{r_1}\;.\label{a55o}
\end{eqnarray}
Both of them correspond to the 2D symmetric semigroup ${\sf S}(r_1,r_3+r_1(r_2-
1))$. In fact, by comparison with (\ref{a51}) we get $k_{10}=k_4$, so we have 
only one new special value $k=k_9$. By comparison the 1st formula in 
(\ref{a55o}) and the 2nd formula in (\ref{a21b}) we obtain,
\begin{eqnarray}
k_9=\left\{\begin{array}{lll}>k_2&\mbox{if}&r_1=2\\=k_2&\mbox{if}&r_1=3\\<k_2
&\mbox{if}&r_1\geq 4\end{array}\right.\;.\label{a55r}
\end{eqnarray} 
%=================================================================
\subsection{Nonsymmetric semigroups ${\sf S}(r_1^2,r_1r_2+r_1^2k,r_3-r_1^2k)$}
\label{s42}
In this section we consider the case of numerical semigroups ${\sf S}(r_1^2,
r_1r_2+r_1^2k,r_3-r_1^2k)$ with nonsymmetric representatives which were not 
observed in sequences with generating triples (\ref{a4}). This case is much 
more difficult to deal with by the reason explained in Remark \ref{rem1}: 
being 3-dimensional,  nonsymmetric semigroup is generated by elements 
satisfying by Lemma \ref{lem1},
\begin{eqnarray}
r_3-r_1^2k\not\in{\sf S}(r_1,r_2+r_1k)\;,\;\;\;\;\;r_3-r_1^2k\leq F(r_1,r_2+
r_1k)\;.\label{a56z}
\end{eqnarray}
The 1st condition in (\ref{a56z}) is necessary and sufficient, however the 2nd 
one is only necessary. Thus, the 2nd condition does not guarantee that the 
chosen $k$ satisfies the 1st one. On the other hand, a straightforward 
application of the latter requirement is hard to perform. 

There exists another problem which makes the construction of nonsymmetric 
semigroups with generating triples (\ref{a50}) not easy. Indeed, summarizing 
(\ref{a31}), (\ref{a22}), (\ref{a54}) and (\ref{a55b}), the set $\Xi\subset{
\mathbb Z}$ of the $k$-values, where nonsymmetric semigroups ${\sf S}(r_1^2,
r_1r_2+r_1^2k,r_3-r_1^2k)$ can be observed, reads
\begin{eqnarray}
\Xi:=\left\{k\;|\;\mu_1<k<\mu_2\right\}\;,\;\;\;\mu_1=\max\left\{
k_1,k_3,k_7\right\},\;\;\mu_2=\min\left\{k_2,k_5,k_8\right\}\;.\label{a56y}
\end{eqnarray}
Thus, if a set $\Xi$ is not empty then every $k_*\in \Xi$ is a candidate to 
make a semigroup with generating triple (\ref{a50}) nonsymmetric. However, what 
is remained open this is a question: does such $k_*$ satisfy the 1st 
requirement in (\ref{a56z}) ?

We can point out the definite values of $k$ associated with nonsymmetric 
semigroups. For example, consider $r_1,r_2,r_3$ providing an integer 
$r_3-r_1^2k$ as a gap of semigroup ${\sf S}(r_1,r_2+r_1k)$,
\begin{eqnarray}
r_3-r_1^2k_*=r_2+r_1k_*+1\;,\;\;\;\;\gcd(r_1,r_2+1)=1\;.\label{a57}
\end{eqnarray}
Equation (\ref{a57}) has the following solution $k_*$ which, by comparison 
with (\ref{a55}), is close to $k_6$,
\begin{eqnarray}
k_*=\frac{r_3-r_2-1}{r_1(r_1+1)}\;,\;\;\;\;\gcd(r_1,r_3)=\gcd(r_1,r_2)=
\gcd(r_1,r_2+1)=1\;.\label{a57a}
\end{eqnarray}
Two last constraints in the right hand side of (\ref{a57a}) forbid $r_1$ be 
divisible by 2. The claim $k_*\in{\mathbb Z}$ requires for $r_2$ and $r_3$ to 
be of distinct parities. It turns out that these properties suffice to give 
rise to infinite family of 2-parametric solutions. Below we give one of them,
\begin{eqnarray}
r_1=2p-1\;,\;\;\;r_2=4p-1\;,\;\;\;r_3=2pk_*(2p-1)+4p\;,\;\;\;p\in
{\mathbb Z}_+\;,\;\;\;p\geq 2\;.\label{a57b}
\end{eqnarray}
In (\ref{a57b}) the value of $k_*$ can be taken on our choice. In Table 1 we 
give a numerical semigroup ${\sf S}(9,21+9k,80-9k)$ which has its nonsymmetric
representatives for $k_*=5,6,7$. In this conjunction, formulas (\ref{a57b}) are 
corresponding to $k_*=6$ and $p=2$ while the other two values of $k_*$ come not 
by (\ref{a57}), but via the other Diophantine equations of similar form.
%=================================================================
\subsection{Concluding Remarks}\label{s43}
In this section we summarize the results on distribution of symmetric and 
nonsymmetric numerical semigroups ${\sf S}(r_1^2,r_1r_2+r_1^2k,r_3-r_1^2k)$ 
governed by one parameter $k$ running throughout the range of its special 
values $k_i$.
\begin{enumerate}
\item In the range $k_4\leq k\leq k_2$ every $k\in {\mathbb Z}$ gives rise to 
the 2D or 3D one parametric numerical semigroups generated by the triple 
(\ref{a50}). 
\item In the range $k_4< k\leq k_3$ all numerical semigroups ${\sf S}(r_1^2,
r_1r_2+r_1^2k,r_3-r_1^2k)$ are symmetric and their minimal generating triple 
cannot be reduced. Their Frobenius numbers coincide with $\Phi(r_1,r_2,r_3)$ 
given by (\ref{a23}).
\item In the range $k_5< k\leq k_2$, $k\neq k_6$, all numerical semigroups 
${\sf S}(r_1^2,r_1r_2+r_1^2k,r_3-r_1^2k)$ are generated by minimal pair $\{
r_1^2,r_3-r_1^2k\}$ and therefore are symmetric. Their Frobenius numbers are 
distinct and differ from $\Phi(r_1,r_2,r_3)$.
\item There exist $k=k_4$ and $k=k_6$ such that the corresponding generating 
sets (\ref{a50}) of semigroups ${\sf S}(r_1^2,r_1r_2+r_1^2k,r_3-r_1^2k)$ are 
reduced up to minimal pairs $\{r_1,r_3-r_1^2k_4\}$ and $\{r_1^2,r_3-r_1^2k_6\}$,
respectively. Their Frobenius numbers coincide with $\Phi(r_1,r_2,r_3)$.
\item In the range $k_1\leq k\leq k_7$ and $k_8\leq k\leq k_2$, $k\neq k_9$, 
all numerical semigroups ${\sf S}(r_1^2,r_1r_2+r_1^2k,r_3-r_1^2k)$ are generated
by minimal pair $\{r_1r_2+r_1^2k,r_3-r_1^2k\}$. Their Frobenius numbers are 
distinct and differ from $\Phi(r_1,r_2,r_3)$.
\item There exists $k=k_9$ such that the corresponding generating set 
(\ref{a50}) is reduced up to minimal pair $\{r_1,r_3+r_1r_2-r_1)$, Its 
Frobenius number coincides with $\Phi(r_1,r_2,r_3)$.
\item In the range $\mu_1<\varkappa<\mu_2$, $\varkappa\in{\mathbb Z}$, and 
$\varkappa\neq k_6$, $\varkappa\neq k_9$, numerical semigroups ${\sf S}(r_1^2,
r_1r_2+r_1^2\varkappa,r_3-r_1^2\varkappa)$ admit their symmetric and 
nonsymmetric representatives, where $\mu_1,\mu_2$ are defined in (\ref{a56y}). 
\end{enumerate}

In Table 1 we present the special values $k_i$ of parameter $k$ for two 
sequences of numerical semigroups discussed in \cite{arn1} and for semigroup 
${\sf S}(9,21+9k,80-9k)$. We give also the Frobenius numbers $F(\left\lfloor 
k_i\right\rfloor)$ associated with these semigroups for $k=\left\lfloor k_i
\right\rfloor$, where $\left\lfloor u\right\rfloor$ denotes the floor function 
of $u$, i.e. $\left\lfloor u\right\rfloor$ gives the largest integer less than 
or equal to $u$.
\begin{center}  
{\bf Table$\;$1.} Semigroups and their Frobenius numbers.
\vspace{.3cm}

\begin{tabular}{|c|c|c|c|c|c|c|c|c|c|c|} \hline\hline
 & $r_1,r_2,r_3$ & $k_1$ & $k_2$ & $k_3$ & $k_4$ & $k_5$ & $k_6$ & $k_7$ & 
$k_8$ & $k_9$\\ \hline\hline
${\sf S}(4,6+4k,87-4k)$ & 2, 3, 87 & -0.5 & 21 & 14.16 & -1 & 15.68 & 14 &-1.24
& 21.49 & 21.25\\\hline
$F(\left\lfloor k_i\right\rfloor)$ & $\Phi$=89 & 89 & 5 & 89 & 89 & 77 & 89 & 
- & 5 & 5 \\\hline\hline
${\sf S}(9,3+9k,85-9k)$ & 3, 1, 85 & 0.33 & 9.11 & 5.66 & 0 & 8.25 & 7 & 
0.21 &9.32 & 9.11 \\ \hline
$F(\left\lfloor k_i\right\rfloor)$ & $\Phi$=167 & 167 & 23 & 167 & 167 & 95 & 
167 & 167 & 23 & 23 \\ \hline\hline
${\sf S}(9,21+9k,80-9k)$ & 3, 7, 80 & -1.66 & 8.55 & 4.53 & -2 & 7.53 & 
6.08 &-2.21 &8.76 & 8.55 \\ \hline
$F(\left\lfloor k_i\right\rfloor)$ & $\Phi$=193 & 193 & 55 & 193 & 193 & 
109 & 121 & - & 55 & 55 \\ \hline\hline
\end{tabular}
\end{center}
In accordance with item 7 of above summary, below we give the values of 
$\varkappa_i$, $\mu_1<\varkappa_i<\mu_2$ and $\varkappa_i\in {\mathbb Z}$ 
together with Frobenius numbers $F(\varkappa_i)$ of corresponding semigroups,
\begin{eqnarray}
{\sf S}(4,6+4\varkappa,87-4\varkappa)&:&\varkappa_1=15,\;F(15)=77\;,\nonumber\\
{\sf S}(9,3+9\varkappa,85-9\varkappa)&:&\varkappa_1=6,\;\varkappa_2=7,\;
F(6)=F(7)=167\;,\;\;\varkappa_3=8,\;F(8)=95\;,\nonumber\\
{\sf S}(9,21+9\varkappa,80-9\varkappa)&:&\varkappa_1=5,\;F(5)=166\;,\;\;
\varkappa_2=6,\;F(6)=121\;,\;\;\varkappa_3=7,\;F(7)=109\;.\nonumber
\end{eqnarray}
The generating sets of all three semigroups are satisfied (\ref{a22c}), i.e. 
$k_1\leq k_3\leq k_2$. Note that in the whole range of varying $k$-parameter, 
$k_1\leq k\leq k_3\leq k_2$ and $k_1\leq k_3\leq k\leq k_2$, including the 
values $\varkappa_i$, both sequences of numerical semigroups ${\sf S}(4,6+4k,
87-4k)$ and ${\sf S}(9,3+9k,85-9k)$ give rise only to symmetric semigroups 
either three-dimensional or two-dimensional. In contrast to them, the numerical 
semigroups ${\sf S}(9,21+9\varkappa_i,80-9\varkappa_i)$ for $\varkappa_i=5,6,7$ 
are three-dimensional and nonsymmetric. 
%=================================================================
\section{Symmetric semigroups in the range $k_1\leq k\leq k_2\leq k_3$}
\label{s5}
In this section we give a detailed analysis on numerical semigroups ${\sf S}(
r_1^2,r_1r_2+r_1^2k,r_3-r_1^2k)$ with parameter $k$ running in the intervals 
(\ref{a22b}) where all semigroups are always symmetric and have the Frobenius 
number $\Phi(r_1,r_2,r_3)$ given by (\ref{a23}). Find $r_1$ such that two 
inequalities in the left and right hand sides in (\ref{a22b}) imposed on 
$r_1r_2+r_3$ become consistent,
\begin{eqnarray}
3+2r_1\leq r_1r_2+r_3\leq\frac{r_1^2+5r_1-3}{r_1-1}\;\;\;\;\rightarrow\;\;\;\;
\frac{r_1(r_1-4)}{r_1-1}\leq 0\;\;\;\;\rightarrow\;\;\;\;2\leq r_1\leq 4\;.
\label{a32}
\end{eqnarray}
Estimate the total number $N$ of such semigroups keeping in mind that according 
to (\ref{a31}) and (\ref{a31a}) $k$ is varying in interval $k_1\leq k\leq  k_2$,
\begin{eqnarray}
N&\leq&\left\lfloor\frac{r_3-3}{r_1^2}\right\rfloor-\left\lfloor\frac{2-r_2}
{r_1}\right\rfloor+1\leq\frac{r_3-3}{r_1^2}-\frac{2-r_2}{r_1}+2=\frac{r_1r_2+
r_3-(2r_1+3)}{r_1^2}+2\nonumber\\
&\leq&\frac1{r_1^2}\left(\frac{r_1^2+5r_1-3}{r_1-1}-2r_1-3\right)+2=
\frac{4-r_1}{r_1(r_1-1)}+2=\frac{(r_1-2)^2}{r_1(r_1-1)}+1< 2\;.\label{a28}
\end{eqnarray}
Thus, a sequence of symmetric numerical semigroups ${\sf S}(r_1^2,r_1r_2+r_1^2k,
r_3-r_1^2k)$ is empty ($N=0$) or contains only one semigroup ($N=1$) for every 
choice of $r_1,r_2,r_3$. In order to find all $k$ providing the case 
(\ref{a22b}) we consider according to (\ref{a32}) all values of $r_1$ 
separately.
\begin{itemize}
\item $r_1=2$, a semigroup ${\sf S}(4,2r_2+4k,r_3-4k)$, $\;2\nmid r_2$, 
$\;2\nmid r_3$.
\begin{eqnarray}
7\leq 2r_2+r_3\leq 11\;,\;\;\;\frac{2-r_2}{2}\leq k\leq\frac{r_3-3}{4}\;.
\label{a33}
\end{eqnarray}
By recasting admitted equations $2r_2+r_3=e_2$ which satisfy the double 
inequality in the left hand side of (\ref{a33}) we have to omit those equations 
when $e_2=0\pmod 2$, otherwise $2\mid r_3$.
\begin{eqnarray}
\left\{\begin{array}{llll}
2r_2+r_3&=&7\;,\;\;&\mbox{has 2 solutions :}\;\{r_2=1,r_3=5\};\;\{r_2=3,r_3=1\}
\\
2r_2+r_3&=&9\;,\;\;&\mbox{has 2 solutions :}\;\{r_2=1,r_3=7\};\;\{r_2=3,r_3=3\}
\\
2r_2+r_3&=&11\;,\;\;&\mbox{has 3 solutions :}\;\{r_2=1,r_3=9\};\;\{r_2=3,r_3=5
\};\;\{r_2=5,r_3=1\}\end{array}\right.\nonumber
\end{eqnarray}
Below we give the corresponding solutions for $k$ and the Frobenius numbers of 
associated numerical semigroups ${\sf S}(4,2r_2+4k,r_3-4k)$ if they exist, i.e.
if $k\in {\mathbb Z}$. We consider two different cases, $2k+r_2\neq 1$ and 
$2k+r_2=1$, or, in other words, when $k$ does satisfy the double inequality in 
the right hand side of (\ref{a33}) and does not satisfy it, respectively.
\begin{eqnarray}
\left.\begin{array}{lllll}
r_2=1,&r_3=7,&k=1,&\;2k+r_2\neq 1,&F(4,6,3)=F(3,4)=5\\
r_2=3,&r_3=3,&k=0,&\;2k+r_2\neq 1,&F(4,6,3)=F(3,4)=5\\
r_2=1,&r_3=9,&k=1,&\;2k+r_2\neq 1,&F(4,6,5)=7\\
r_2=3,&r_3=5,&k=0,&\;2k+r_2\neq 1,&F(4,6,5)=7\\
r_2=5,&r_3=1,&k=-1,&\;2k+r_2\neq 1,&F(4,6,5)=7\end{array}\right.\label{a40}
\end{eqnarray}
\begin{eqnarray}
\left.\begin{array}{lllll}
r_2=1,&r_3=5,&k=0,&\;2k+r_2=1,&F(4,2,5)=F(2,5)=3\\
r_2=3,&r_3=1,&k=-1,&\;2k+r_2=1,&F(4,2,5)=F(2,5)=3\\
r_2=1,&r_3=7,&k=0,&\;2k+r_2=1,&F(4,2,7)=F(2,7)=5\\
r_2=3,&r_3=3,&k=-1,&\;2k+r_2=1,&F(4,2,7)=F(2,7)=5\\
r_2=1,&r_3=9,&k=0,&\;2k+r_2=1,&F(4,2,9)=F(2,9)=7\\
r_2=3,&r_3=5,&k=-1,&\;2k+r_2=1,&F(4,2,9)=F(2,9)=7\\
r_2=5,&r_3=1,&k=-2,&\;2k+r_2=1,&F(4,2,9)=F(2,9)=7\end{array}\right.
\label{a40a}
\end{eqnarray}
\end{itemize}
\begin{itemize}
\item $r_1=3$, a semigroup ${\sf S}(9,3r_2+9k,r_3-9k)$, $\;3\nmid r_2$, 
$\;3\nmid r_3$,
\begin{eqnarray}
9\leq 3r_2+r_3\leq \frac{21}{2}\;,\;\;\;\frac{2-r_2}{3}\leq k\leq\frac{r_3-3}
{9}\;.\label{a34}
\end{eqnarray}
Omit equations $3r_2+r_3=e_3$ such that $e_3=0\pmod 3$, otherwise $3\mid r_3$. 
Thus, we have,
\begin{eqnarray}
3r_2+r_3=10\;,\;\;&\mbox{has 2 solutions :}\;\{r_2=1,r_3=7\};\;
\{r_2=2,r_3=4\}\;.\label{a37}
\end{eqnarray}
Similarly to the previous case we give the corresponding solutions for $k$, 
$k\in {\mathbb Z}$, and the Frobenius numbers of associated numerical 
semigroups ${\sf S}(9,3r_2+9k,r_3-9k)$ in two different cases, $3k+r_2\neq 1$ 
and $3k+r_2=1$. 
\begin{eqnarray}
\left.\begin{array}{lllll}
r_2=2,&r_3=4,&k=0,&\;3k+r_2\neq 1,&F(9,6,4)=11\;,\\\label{a38}
r_2=1,&r_3=7,&k=0,&\;3k+r_2=1,&F(9,3,7)=F(3,7)=11.\end{array}\right.\label{a38}
\end{eqnarray}
\end{itemize}
\begin{itemize}
\item $r_1=4$, a semigroup ${\sf S}(16,4r_2+16k,r_3-16k)$, $\;2\nmid r_2$, 
$\;2\nmid r_3$,
\begin{eqnarray}
4r_2+r_3=11\;,\;\;\;\frac{2-r_2}{4}\leq k\leq\frac{r_3-3}{16}\;.\label{a35}
\end{eqnarray}
It turns out that in the case $4k+r_2\neq 1$ an equation (\ref{a35}) has not 
an integer solution in $k$. Thus, the only numerical semigroup ${\sf S}
(16,4r_2+16k,r_3-16k)$ with corresponding Frobenius number reads,
\begin{eqnarray}
\left.\begin{array}{lllll}
r_2=1,\;\;r_3=7,\;\;k=0,\;\;4k+r_2=1,\;\;F(16,4,7)=F(4,7)=17\;.
\end{array}\right.
\label{a36}
\end{eqnarray}
\end{itemize}
%=================================================================
\section{Symmetric semigroups ${\sf S}(r_1^2,r_3-r_1^2k)$ and enumeration of\\
integer points in plane curve}\label{s6}
In section \ref{s41} we have observed a phenomenon of reduction of a number of
three minimal generators $\{r_1^2,r_1r_2+r_1^2k,r_3-r_1^2k\}$ up to two due to 
the linear dependence between the 1st and the 2nd generators if $r_1k+r_2=1$. 
In this conjunction ask about the other ways of similar reduction of semigroup's
dimension when the linear dependence arises in the rest of two pairs of 
generators separately, namely, between the 1st and the 3rd generators or 
between the 2nd and the 3rd generators.

Regarding the first pair of the 1st and the 3rd generators, it can be proven 
that by assumptions $\gcd(r_1,r_3)=1$ and $r_3\in {\mathbb Z}^+$ such linear 
dependence could not happen. Indeed, if such dependence holds, $r_3-r_1^2k=c
r_1^2$, $c\in{\mathbb Z}^+$, then $r_3$ is divisible by $r_1^2$ or vanishes 
that contradicts the above assumptions.

Regarding the second pair of the 2nd and the 3rd generators, their linear 
dependence
\begin{eqnarray}
r_3-r_1^2k=f\cdot (r_1r_2+r_1^2k)\;,\;\;\;f\in{\mathbb Z}^+\;,\label{a43a}
\end{eqnarray}
could not happen since it also contradicts the assumptions $\gcd(r_1,r_3)=1$ 
and $r_3\in {\mathbb Z}^+$.

Thus, consider the following linear dependence,
\begin{eqnarray}
r_1r_2+r_1^2k=g\cdot (r_3-r_1^2k)\;,\;\;\;g\in{\mathbb Z}^+\;.\label{a43}
\end{eqnarray}
The quadratic Diophantine equation (\ref{a43}) describes an algebraic curve of 
degree 2 in the $k-g$ plane. The number of points with integer coordinate, 
$k\in{\mathbb Z}$ and $g\in{\mathbb Z}^+$, of this curve coincides with a 
number of solutions of the Diophantine equation (\ref{a43}). It can be solved 
completely by reduction it to the Pell equation and further calculation of 
continued fractions \cite{wh12}.

In this section we give necessary conditions to have the integer solutions, 
$k\in{\mathbb Z}$ and $g\in{\mathbb Z}^+$, of equation (\ref{a43}) and present 
two examples associated with Arnold's experiments showing how these 
requirements help to find all triples with linear dependence between the 2nd 
and the 3rd generators.

First, note that $g$ is divisible by $r_1$ that follows by (\ref{a43}) and 
assumption $\gcd(r_1,r_3)=1$. Denote $X=g+1$, $Y=k$ and rewrite equation 
(\ref{a43}) as follows,
\begin{eqnarray}
r_1^2Y-r_3+\frac{r_1r_2+r_3}{X}=0\;,\;\;\;Y\in{\mathbb Z},\;\;X\in{\mathbb 
Z}^+,\;\;X\geq 2\;.\label{a44}
\end{eqnarray}
The Diophantine equation (\ref{a44}) is solvable iff $X$ takes its value 
among divisors of $r_1r_2+r_3$. Hence the next Lemma follows.
\begin{lemma}\label{lem2}
Let a numerical semigroup ${\sf S}(r_1^2,r_1r_2+r_1^2k,r_3-r_1^2k)$ be given 
such that $k\in{\mathbb Z}$, $r_1,r_2,r_3\in {\mathbb Z}^+$ and $r_1\geq 2$, 
$\gcd(r_1,r_2)=\gcd(r_1,r_3)=1$. If the linear dependence (\ref{a43}) holds 
for $g=g_*$ and $k=k_*$ then a semigroup is isomorphic to the 2D symmetric 
semigroup ${\sf S}(r_1^2,r_3-r_1^2k_*)$, and $g_*$ is divisible by $r_1$, 
and $r_1r_2+r_3$ is divisible by $g_*+1$. 
\end{lemma}

Denote by $Q(r_1,r_2,r_3)$ the total number of solutions of equation 
(\ref{a44}) and by $\sigma_0(n)$ the number of positive divisors $\delta_i(n)$ 
of integer $n$, where $i=1,\ldots,\sigma_0(n)$. First, by $g=X-1\geq 1$ we 
have to exclude the minimal divisor $\delta_{min}(r_1r_2+r_3)=1$ from possible 
solutions $X$ of (\ref{a44}). Next, let, by way of contradiction, the maximal 
divisor $\delta_{max}(r_1r_2+r_3)=r_1r_2+r_3$ coincides with one of solutions 
$X$. Then substituting it into (\ref{a44}) we get the final triple: $\{r_1^2,
r_1r_2+r_3-1,1\}$. The occurrence of unity in the minimal generating set, 
$1\in {\bf d}^3$, makes the associated numerical semigroup ${\sf S}\left(
{\bf d}^3\right)$ free of gaps and equivalent to the whole set of nonnegative 
integers, ${\sf S}\left({\bf d}^3\right)\equiv{\mathbb Z}^+\cup\{0\}$. For 
such semigroup the Frobenius number does not exist.

Thus, there exist $Q(r_1,r_2,r_3)$ different sporadic values $k=k_*$ which 
suffice to reduce the dimension of numerical semigroups ${\sf S}(r_1^2,r_1r_2+
r_1^2k,r_3-r_1^2k)$ up to 2 and induce a bijective correspondence between 
symmetric semigroups, ${\sf S}(r_1^2,r_1r_2+r_1^2k_*,r_3-r_1^2k_*)
\leftrightarrow {\sf S}(r_1^2,r_3-r_1^2k_*)$, with nonempty sets of gaps. 
Keeping in mind both values $\delta_{min}(r_1r_2+r_3)$ and $\delta_{max}(r_1r_2
+r_3)$, we can give the lower and upper bounds for $Q(r_1,r_2,r_3)$,
\begin{eqnarray}
0\leq Q(r_1,r_2,r_3)\leq\sigma_0(r_1r_2+r_3)-2\;.\label{a45}
\end{eqnarray}
By (\ref{a45}) and Lemma \ref{lem2} we come to the other Corollaries related 
to the cases when $Q(r_1,r_2,r_3)=0$.
\begin{corollary}\label{cor1}
Let a numerical semigroup ${\sf S}(r_1^2,r_1r_2+r_1^2k,r_3-r_1^2k)$ be given as 
in Lemma \ref{lem2} and $r_1r_2+r_3$ is a prime number. Then one cannot choose 
$g=g_*$ and $k=k_*$ such that the linear dependence (\ref{a43}) does define the 
2D semigroup ${\sf S}(r_1^2,r_3-r_1^2k_*)$ with nonempty set of gaps.
\end{corollary}
\begin{corollary}\label{cor2}
Let a numerical semigroup ${\sf S}(r_1^2,r_1r_2+r_1^2k,r_3-r_1^2k)$ be given 
as in Lemma \ref{lem2} and $r_1r_2+r_3=p^2$ where $p$ is a prime number. 
Then $Q(r_1,r_2,r_3)=0$ if $p-1$ is not divisible by $r_1$.
\end{corollary}
{\sf Proof} $\;\;\;$By (\ref{a45}) there is only one candidate for solutions 
of the Diophantine equation (\ref{a44}) and by Corollary \ref{cor2} this is 
$X=p$. Substituting it into (\ref{a44}) we get
\begin{eqnarray}
r_1^2Y=p(p-1)-r_1r_2\;,\;\;\;\;\;\mbox{or}\;\;\;\;\;r_1(r_2+r_1Y)=p(p-1)\;.
\label{a47}
\end{eqnarray}
Note that $\gcd(r_1,p)=1$, or, keeping in mind that $p$ is a prime number, 
this is equivalent that $r_1$ is not divisible by $p$. Indeed, let, by way of 
contradiction, $r_1=v\cdot p$, $v\in{\mathbb Z}^+$. Then $r_3$ is also 
divisible by $p$ since, by $r_1r_2+r_3=p^2$, we have $r_3=p(p-vr_2)$. However, 
this contradicts the assumptions $\gcd(r_1,r_3)=1$ and $r_3\in {\mathbb Z}^+$. 

If the Diophantine equation (\ref{a47}) is solvable then $p-1$ is necessarily 
divisible by $r_1$. Thus, if $p-1$ is not divisible by $r_1$ then equation 
(\ref{a47}) is not solvable, i.e. $Q(r_1,r_2,r_3)=0$.$\;\;\;\;\;\;\Box$

Corollary \ref{cor2} gives the necessary but not sufficient conditions for 
equation (\ref{a47}) to be solvable. Indeed, let $p-1=u\cdot r_1$, $u\in
{\mathbb Z}^+$. Substituting it into (\ref{a47}) we get $r_1Y=u\cdot p-r_2$. 
Thus, the solvability of the last Diophantine equation in $Y$ presumes an 
additional divisibility relation.

In the following Examples \ref{ex1} and \ref{ex2} we present the phenomenon of 
reduction of the 3D semigroup's dimension up to 2 in two different sequences 
of semigroups generated by triples (\ref{a4}) and discussed in \cite{arn1}. In 
both Examples we have underlined those divisors $\delta_i$ of $r_1r_2+r_3$ 
which give rise to sporadic 2D semigroups with corresponding $k_*$ and $g_*$. 
\begin{example}$\{d_1,d_2,d_3\}=\{4,6+4k,87-4k\}$, $\{r_1,r_2,r_3\}=
\{2,3,87\}$\label{ex1}
\begin{eqnarray}
&&r_1r_2+r_3=93\;,\;\;\;\sigma_0(93)=4\;,\;\;\;\delta_i(93)=1,\underline{3},
\underline{31},93\;,\;\;\;Q(2,3,87)=2\;,\nonumber\\
&&r_1\mid\delta_i(93)-1\;:\;\;2\mid 2,\;\;2\mid 30,\;\;2\mid 92\;,\nonumber\\
&&(k_{1*},g_{1*})=(14,2),\;\;{\sf S}\left({\bf d}^3_1\right)={\sf S}(4,62,31),
\;\;F(4,62,31)=F(4,31)=89\;,\nonumber\\
&&(k_{2*},g_{2*})=(21,30),\;\;{\sf S}\left({\bf d}^3_2\right)={\sf S}(4,90,3),
\;\;F(4,90,3)=F(4,3)=5\;.\nonumber
\end{eqnarray}
\end{example}
\begin{example}$\{d_1,d_2,d_3\}=\{9,3+9k,85-9k\}$, $\{r_1,r_2,r_3\}=
\{3,1,85\}$\label{ex2}
\begin{eqnarray}
&&r_1r_2+r_3=88\;,\;\;\;\sigma_0(88)=8\;,\;\;\;\delta_i(88)=1,2,
\underline{4},8,11,\underline{22},44,88\;,\;\;\;Q(3,1,85)=2\;,\nonumber\\
&&r_1\mid\delta_i(88)-1:\;\;3\mid 3,\;\;3\mid 21,\;\;3\mid 87\;,\nonumber\\
&&r_1\nmid\delta_i(88)-1:\;\;3\nmid 1,\;\;3\nmid 7,\;\;3\nmid 10,\;\;3\nmid 43
\;,\nonumber\\
&&(k_{1*},g_{1*})=(7,3),\;\;{\sf S}\left({\bf d}^3_1\right)={\sf S}(9,66,22),
\;\;F(9,66,22)=F(9,22)=167\;,\nonumber\\
&&(k_{2*},g_{2*})=(9,21),\;\;{\sf S}\left({\bf d}^3_2\right)={\sf S}(9,84,4),
\;\;F(9,84,4)=F(9,4)=23\;.\nonumber
\end{eqnarray}
\end{example}
Regarding the 3rd semigroup ${\sf S}(9,21+9k,80-9k)$, where nonsymmetric
representatives are admitted (see section \ref{s42}), we have $r_1r_2+r_3=101$,
$\sigma_0(101)=2$ that by Corollary \ref{cor1} results in $Q(3,7,80)=0$, i.e.
the Diophantine equation (\ref{a44}) has no solutions.
%=================================================================
\section{Acknowledgement}
The research was partly supported by the Kamea Fellowship.
%=================================================================

\end{document}